\documentclass[12pt]{article}
\usepackage{latexsym,amsfonts,dif}
\usepackage[round]{natbib}

\newcommand{\sect}[1]{\section{#1} \setcounter{equation}{0}
                      \setcounter{table}{0} \setcounter{figure}{0}}

\newtheorem{defi}{Definition}[section]
\newtheorem{ex}[defi]{Example}
\newtheorem{con}[defi]{Condition}
\newtheorem{lem}[defi]{Lemma}
\newtheorem{thm}[defi]{Theorem}

\newenvironment{pr}{\ni {\bf Proof}.~~}{}  

\newtheorem{assu}[defi]{Assumption}

\title{A review of asymptotic theory of estimating functions}

\author{Jean Jacod \\ 
{\small Institut de Math\'ematiques de Jussieu and \vspace{-1.5mm}} \\
{\small Universit\'e P.\ et M.\ Curie (Paris-6) \vspace{-1.5mm}} \\
{\small CNRS UMR 7586 \vspace{-1.5mm}} \\
{\small 75252 Paris C\'edex 05 \vspace{-1.5mm}} \\
{\small France}
\and
Michael S\o rensen \\ 
{\small Dept.\ of Mathematical Sciences \vspace{-1.5mm}} \\
{\small University of Copenhagen\vspace{-1.5mm}} \\
{\small Universitetsparken 5\vspace{-1.5mm}} \\
{\small DK-2100 Copenhagen {\O}\vspace{-1.5mm}} \\
{\small Denmark}}

\date{}

\begin{document}

\maketitle

\begin{abstract}
Asymptotic statistical theory for estimating functions is reviewed in
a generality suitable for stochastic processes. Conditions concerning
existence of a consistent estimator, uniqueness, rate of convergence,
and the asymptotic distribution are treated separately. Our conditions
are not minimal, but can be verified for many interesting stochastic
process models. Several examples illustrate the wide applicability of
the theory and why the generality is needed. \vspace{1mm}
\\ 
{\bf Key words:} Asymptotic statistics, diffusion processes, ergodic processes, high
frequency asymptotics, limit theory, longitudinal data, Markov
process, misspecified model.

\end{abstract}


\sect{Introduction}

Estimating functions provide a general framework for finding estimators
and studying their properties in many different kinds of statistical
models and asymptotic scenarios. 
We present the main results of the asymptotic theory of estimating
functions in a generality that is suitable for statistical inference
for stochastic processes. Typically, we have observations of
$n$ random variables, $X_1, \ldots, X_n$, whose joint law depends on a
parameter $\te$. Also, for instance, the sampling frequency or
the diffusion coefficient of a diffusion process may depend on $n$. 
An estimating function is a function
of the data as well as the parameter, $G_n(\theta) = G_n(\theta, X_1,
\ldots, X_n)$, and estimators are obtained by solving the estimating
equation $G_n(\theta) =0$. We present simple conditions for existence,
uniqueness, consistency, rate of convergence and asymptotic
distribution of such estimator as $n$ tends to infinity. Our conditions are
not minimal, but can be verified for many interesting stochastic
process models. This is illustrated by several examples, which also
demonstrate why the generality of our conditions is needed.

The theory covers
consistent estimators obtained by maximising (or minimising), with
respect to $\theta$, 
a function $H_n(\theta, X_1, \ldots, X_n)$ of the data and the
parameter. If $H_n$ is continuously differentiable, and if the true
parameter value belongs to the interior of the parameter space, then
eventually the estimator will solve the estimating equation
$\partial_\theta H_n(\theta)=0$. This includes, of course, maximum
likelihood estimation. It is often the case for stochastic process
models that the likelihood function is not explicitly available, is
not tractable, or requires heavy computations. In such cases other
estimating functions provides a tractable alternative.

The idea of using estimating equations is an old one and goes back at
least to Karl Pearson's introduction of the method of moments. The
term estimating function may have been coined by \cite{kimball}. 
In the statistics literature, the modern theory of optimal estimating
functions dates back to the papers by \cite{godambe60} and
\cite{durbin}, however the basic idea was in a sense already used in
\cite{fisher}. The theory was extended to stochastic processes 
by \cite{godambe85}, \cite{godambeheyde}, and several others; 
see the references in \cite{heyde97}. In the econometrics literature
the foundation was laid by \cite{hansen82}. This paper treats a general
class of models including time series models. Asymptotic theory was
developed in, e.g., \cite{hansen82,hansen2}, \cite{chamberlain87}, and  
\cite{kuersteiner}. The theory was extended by several authors, see
the discussion and the discussion and references in \cite{hall}.
A discussion of links between the econometrics and statistics
literature can be found in \cite{hansen2001}.

The classical approach to asymptotic statistical theory for estimating
functions is based on the seminal work of \cite{cramer}. To prove
asymptotic existence of an estimator, one approach, originally due to
\cite{aitchison}, is based on Brouwer's fixed-point theorem. This
idea is used in, e.g., \cite{sweeting}, \cite{oebnms}, and 
\cite{ms00a}. In the present paper, we follow instead an approach based on
the fixed-point theorem for contractions, which was also used by, e.g.,
\cite{yuanjennrich}. The latter approach gives relatively simple
conditions and a straightforward identifyability condition that
implies a uniqueness result for estimators. A comprehensive
introduction to asymptotic statistics can be found in \cite{vaart}.

The paper is organised as follows. In Section 2 we present the
general asymptotic theory, while theory for ergodic processes is
presented in Section 3. Longitudinal data are briefly considered in
Section 4, and high frequency data, both with a fixed and an infinite
time horizon, are discussed in Section 5. All proofs are collected in
Section 6.

\sect{Asymptotics for estimating functions}
\label{asyest}

The general set-up is as follows. We have a measurable space
$(\Om,\f)$ with a probability measure $\PP$, the {\em
true measure}, and a family $(\PP\itt)_{\te\in\Te}$ of probability
measures indexed by $\Te\subseteq\R^p$, the {\em statistical
  model}. Often it is assumed that there exists a $\theta_0 \in
\Theta$ (called the true parameter value) such that $\PP=\PP\ittt$,
but this need not be the case. 

At stage $n$, we have a set of observations which generates a $\PP$--complete
$\si$--field $\f_n$ (that is, $\f_n$ is the $\PP$--completion of the
$\si$--field $\f^0_n$ generated by those observations.) Taking the
completion simplifies the mathematical formulation of the results below,
and it has no negative consequences from a practical viewpoint since all statements
about estimators are always up to a $\PP$-null set. To be
consistent, we also suppose that $\f$ is $\PP$--complete. If we observe a
continuous time stochastic process
$X_t$, then $\f_n$ could, for instance, be the (complete) $\si$-field
generated by the variables $X_s$ for $s\in[0,t_n]$ for some increasing
sequence $t_n$, or by the variables $X_{i\De_n}$ for $i=0,\ldots,n$
and some $\De_n>0$. 

An {\em estimating function} at stage $n$ is a function
$(\te,\om)\mapsto G_n(\te,\om)$ that takes its values in $\R^p$ 
and depends on the statistical parameter
$\theta \in \Theta$ and on the observation at this stage, that
is $(\te,\om)\mapsto G_n(\te,\om)$ is measurable w.r.t.\ the product
of the Borel $\si$-field of $\Te$ with $\f_n$. For convenience, we
usually suppress the argument $\om$ in the notation. We get an
estimator $\tten$ by solving the estimating equation
\bee
\label{esteq3.1}
G_n (\theta) = 0.
\eee
For any given $\om$, this equation may have a unique solution, 
several solutions, or none, so we have to be a bit careful. Therefore
we give the following formal definition, where $\de$ denotes a
``special'' point, which we take to be outside $\Te$ because this is
most convenient, and $\Te_\de= \Te\cup\{\de\}$.

\begin{defi}\label{EST-FU} a) The {\em domain of definition} of
  $G_n$-estimators (for a given $n$) is the set $D_n$ of all $\om$ for which
  $G_n(\te,\om)=0$ for at least one $\te\in\Te$.

b) A {\em $G_n$-estimator}, $\tten$, is any $\f_n$--measurable map from $\Om$
into $\Te_\de$, such that $\tten\in\Te$ and $G_n(\tten)=0$ on the
set $D_n$ and $\tten=\de$ on its complement $(D_n)^c$.
\end{defi}

Because $\f_n$ is $\PP$-complete, the measurable
selection theorem implies that $A_n\in\f_n$ and that a $G_n$--estimator
always exists.

In the rest of this paper, $\m_p$ denotes the set of all $p\times p$ matrices, and
$\m_p^{\rm inv}$ the subset of all $A\in\m_p$ that are invertible (non--singular).
A $G_n$--estimator is also a $B_nG_n$--estimator for any
$B_n\in\m_p^{\rm inv}$, where $B_n$ can depend on the data as well as the
parameter. We call the family of estimating functions of the form $B_nG_n$
{\em versions} of $G_n$. Since $B_n$ depends on $n$, the assumptions
about the asymptotic behavior of $G_n$ made in the following can obviously
fail for some versions of $G_n$, but it is sufficient that a version exists, for
which the conditions are satisfied. All limits below are taken
as $n\rightarrow\infty$.

\subsection{Existence and uniqueness of a consistent estimator}
\label{existence}

The following condition ensures that
for $n$ large enough the estimating equation has a solution that converges
to a particular parameter value $\bar{\theta}$. When the statistical model
contains the true model, the estimating function should preferably be chosen
such that $\bar{\theta} = \theta_0$. It is, however, useful to include
the more general case in the theory. A couple of exemples below will
illustrate why.

\begin{con}
\label{conasy1}
There exists a parameter value $\tth\in \mbox{\em int} \,\Te$ (the interior
of $\Te$), a neighbourhood $M$ of $\tth$, and a (possibly random)
$\R^p$--valued function $G$ on $M$, such that the following holds:

{\em (i)} $G_n (\bar{\theta})~\stackrel{\PP}{ \rightarrow}~
0$ (convergence in probability under the true measure $\PP$) and $G(\tth)=0$.

{\em (ii)} For $\PP$-almost all $\om$, $G(\cdot,\om)$ and all
$G_n(\cdot,\om)$ are $C^1$ (i.e., continuously differentiable) on $M$, and
\bee
\label{unifcon}
\sup_{\te\in M}\|\,\partial_\te G_n (\te) -\partial_\te G(\te) \, \|
~\stackrel{\PP}{ \rightarrow}~ 0.
\eee

{\em (iii)} The matrix $\partial_\te G(\tth)$ is
non-singular with $\PP$--probability one.
\end{con}

In (\ref{unifcon}) $\|\cdot\|$ could be any norm on $\m_p$, since
any two norms on a finite-dimensional space are equivalent. However, in the
following we will use the operator norm
$$\|A\|^2=\max\{|\la_j|:\la_j\,\mbox{ is an eigenvalue of } A^* A \},$$
or equivalently $\|A\|=\sup_{|x|=1}|A x|.$ Here and below $|\cdot|$ denotes
the Euclidian norm, and $A^*$ is the transpose of $A$. The vector
$\bar{\te}$ and the function $G(\te)$ depend on
the true probability measure $\PP$. 
Note that when $M$ is a bounded set, (\ref{unifcon}) implies
\bee
\label{unifcon'}
\sup_{\theta \in M} |G_n (\theta)- G (\theta) | \stackrel{\PP}{\rightarrow}0.
\eee

It is not a restriction to assume in Condition \ref{conasy1} that $M$
is convex, or even an open or closed ball 
centered at $\bar{\te}$, since $M$ always contains such a ball.
Since $G_n$ and $G$ are $C^1$ on $M$ outside a
$\PP$-null set, whereas $\f$ is $\PP$-complete, the left side of
(\ref{unifcon}) is $\f$--measurable. Moreover, it is not necessary to
assume the existence of the function $G$ in 
the condition, where $\partial_\te G(\te)$ can be replaced
by a function $W(\te) \in \m_p$ such that $W(\bar \te) \in \m_p^{\rm inv}$. 
Then the existence of a $G$ such that $\partial_\te G(\te) = W(\te)$
and $G(\bar \te)=0$ follows.

At this point it is not obvious why we do not simply take $\bar{\te}$ to
be equal to the true parameter value $\theta_0$. Before presenting the
asymptotic theory, let us therefore present two simple examples showing
that there are situations where quite naturally $\bar{\te} \neq \theta_0$.
Note that the theory below would in no way be
simplified by assuming that  $\bar{\te} = \theta_0$.

\begin{ex}
\label{exasy1}
{\em Suppose we model observations $X_0, X_1, \ldots , X_n$ by
the autoregression of order one
\bee\label{AR1}
X_i = \theta X_{i-1} + \epsilon_i,
\eee
where the $\epsilon_i$'s are i.i.d.\ random
variables with mean zero and finite variance, and $\te \in (-1,1)$ so that
$X$ is ergodic. It is natural to estimate $\theta$ by minimizing
$$H_n(\theta) = n^{-1} \sum_{i=1}^n (X_i - \theta X_{i-1})^2,$$
which is minus the logarithm of a Gaussion pseudo-likelihood.
This least squares estimator can be found by solving the estimating equation
$G_n(\te)=0$, where
\[
G_n(\theta) = n^{-1 }\sum_{i=1}^n X_{i-1} (X_i - \theta X_{i-1}).
\]
If our observations are in fact generated by
(\ref{AR1}) with $\te=\te_0$, then $G_n(\te_0)=n^{-1}\sn X_{i-1}
\epsilon_i$ $\stackrel{\PP}{ \rightarrow} 0$ by the law of large numbers
for martingales. It is not difficult to see that $\te_0$ is the only
parameter value for which $G_n(\te)\stackrel{\PP}{ \rightarrow} 0$,
so necessarily $\bar{\te}=\te_0$.

Now assume that our data are actually observations from an autoregression
of order two, i.e.\ that
\[
X_i = \theta_1 X_{i-1} + \theta_2 X_{i-2} + \epsilon_i,
\]
where the $\epsilon_i$'s are as before, and where $\theta_1$ and $\theta_2$
are such that the observed process $X$ is ergodic. Thus if $\theta_2 \neq
0$, our statistical model is misspecified. When this is the case,
\[
G_n(\theta) = (\theta_1 - \theta) n^{-1 }\sum_{i=1}^n X_{i-1}^2 +
\theta_2  n^{-1}\sum_{i=1}^n X_{i-1}X_{i-2} +
n^{-1 }\sum_{i=1}^n X_{i-1} \epsilon_i,
\]
implying that
\[
G_n(\theta) \stackrel{\PP}{ \rightarrow} (\theta_1 - \theta) \sigma^2
+ \theta_2 \nu,
\]
where $\si^2$ and $\nu$ are the expectations of $X_i^2$ and $X_{i-1}X_i$
under the stationary distribution for $X$. We see that necessarily
$\bar{\te}=\te_1+\te_2\nu/\si^2$. The
least squares estimator converges to $\bar{\te}$ as $n \rightarrow \infty$
when the observed process $X$ is an autoregression of order two.}\qed
\end{ex}

\begin{ex}
\label{exasy2}
{\em
As another closely related example, consider the Ornstein-Uhlenbeck
process given by 
\bee
\label{asyOU}
dX_t = -\te X_t dt + d W_t, \ \ \ \te > 0.
\eee
Suppose we have observations $X_0, X_\Delta, \ldots , X_{n \Delta}$
from (\ref{asyOU}) with $\theta = \theta_0$. The parameter $\te$ is often
estimated by the least squares estimator obtained by minimizing the
function
\[
H_n(\theta) = n^{-1} \sum_{i=1}^n (X_{i\Delta} - (1-\theta \Delta)
X_{(i-1)\Delta})^2
\]
This is motivated by an Euler discretization and is
the same as using the estimating function 
\[
G_n(\te) = n^{-1} \sum_{i=1}^n X_{(i-1)\Delta}[X_{i\Delta} - (1-\theta \Delta)
X_{(i-1)\Delta}].
\]
Since $X_{i\Delta} = e^{-\theta_0\Delta} X_{(i-1)\Delta} + \epsilon_i,$
where the $\epsilon_i$'s satisfy the assumptions in Example \ref{exasy1},
it follows from this example that $(1-\bar{\te}\Delta) = e^{-\theta_0\Delta}$
or $\bar{\te} = (1 - e^{-\theta_0\Delta})/\Delta$.
The mean squares estimator $\hat{\te}_n$ obtained by solving
$G_n(\theta)=0$ converges to $\bar{\te}$ as $n \rightarrow \infty$. When
$\theta_0\Delta$ is small, $\bar{\te}$ is close to $\te_0$, but if
$\theta_0\Delta$ is large, $\bar{\te}$ will be very far from $\te_0$.
}
\qed

\end{ex}

Condition \ref{conasy1} implies that $\bar{\te}$ is a.s.\ an isolated
zero of $G$. In view of (\ref{unifcon'}), it is intuitively
clear that $G_n $ must have a unique zero near $\bar{\theta}$ when
$n$ is sufficiently large. That is the content of the following theorem.
When $\theta$ is one-dimensional, it is easy to turn the intuitive
argument into a rigorous proof. When $p > 1$, this is most easily done using
a fixed point theorem as we shall do in the proof of Theorem \ref{thmasy1}.

A sequence of estimators that converges to a parameter value different
from the true value is usually not called consistent, but to facilitate the
discussion, we call a sequence $(\tten)$ of estimators {\em weakly
  $\bar{\te}$-consistent} if $\tten~ \stackrel{\PP}{\longrightarrow}~
\bar{\theta}$ as $n\to\infty$. If the convergence is
almost sure w.r.t.\ $\PP$, the estimator is called {\em strongly
  $\bar{\te}$-consistent}. 

Since $\tten$ takes its values in the ``extended'' space $\Te_\de$,
one has to define a
metric on this set in order that the previous convergence makes
sense. In $\Te$ it is of course the restriction of the usual metric on
$\R^p$; and we will say that the distance between the extra point $\de$
and any point in $\Te$ equals $1$. Therefore weak $\bar{\te}$-consistency 
implies in particular that $\PP(\tten=\de)\to0$, or otherwise formulated:
$\tten$ solves the estimating equation (\ref{esteq3.1}) with a probability
tending to one.

\begin{thm}
\label{thmasy1}
Under Condition \ref{conasy1} we can find a sequence $(\tten)$ of
$G_n$-estimators which is weakly $\bar{\te}$-consistent.
Moreover this sequence is {\em eventually unique}, that is if
$(\tten')$ is any other weakly $\bar{\te}$--consistent sequence of
$G_n$--estimators, then $\PP(\tten\neq\tten')\to0$ as $n\to\infty$.
\end{thm}

This theorem is proved in Section \ref{proofs}, and a similar proof gives
the following result.

\begin{thm}
\label{thmasy2}
Suppose Condition \ref{conasy1} holds with almost sure convergence
instead of convergence in probability. Then a sequence $(\tten)$ of
$G_n$--estimators exists which is {\em strongly $\bar{\te}$--consistent}.
\end{thm}

It is important to observe that so far we have made no assumption about
identifiability. Consequently, although there is a sequence of
$G_n$--estimators that is (weakly or strongly) $\bar{\te}$--consistent,
there might be other sequences that are not. For example, if $G$
vanishes for another value $\bar{\te}'$ and $\partial_\te
G(\bar{\te}') \in \m_p^{\rm inv}$,  then there is another sequence of
$G_n$--estimators that  
converges to $\bar{\te'}$. An estimating function with such a property is
obviously not of much practical use. Also, the reader can observe that
the construction of $\tten$ in the proof of the previous theorems made
use of the value $\bar{\te}$. Thus this result is a mathematical
existence result, and the method can obviously not be used in
statistical practice to choose a good $G_n$--estimator among several
solutions to $G_n(\te)=0$ since $\bar{\te}$ is unknown. 

Thus the previous two theorems are powerful technical tools, but to
obtain practically useful results, global properties of the estimating
function are needed. Essentially, $\tth$ must be the only root to
$G(\te)=0$. Global uniqueness results can be obtained when $\Theta$ is
compact and when $\Te$ is homeomorphic to $\R^p$ (and hence is open).
These two cases cover most practical situations. Recall that $\Te$
is homeomorphic to $\R^d$ if there is a bijective bicontinuous mapping
$\psi:\,\Te\rightarrow\R^p$. 

When $\Te$ is homeomorphic to $\R^p$, to obtain a global
uniqueness result, we must restrict attention to a
class of estimators that are prevented from going to the boundary of
$\Te$. We define such a class as follows. Pick an arbitrary point 
$\rho\in\Te$, and define a {\em $\rho$--centered $G_n$--estimator} as
any $\f_n$--measurable and $\Te_\de$-valued variable $\tten^\rho$ such
that
\bee\label{CC1}
\tten^\rho=\left\{\begin{array}{ll}
\mbox{argmin}\,\big(d(\te,\rho):\,\te\in\Te,\,G_n(\te)=0\big)\qquad
& \mbox{on the set } D_n \\
\de & \mbox{on } (D_n)^c,
\end{array}\right.
\eee
where $d(\te,\te')=|\psi(\te)-\psi(\te')|$ ($\psi$ is the bijection
$\Te\rightarrow\R^p$), and $\delta$ and $D_n$ are 
as in Definition \ref{EST-FU}. By the measurable selection theorem and
the continuity of $G_n$, such estimators always exist.

\begin{thm}
\label{thmasy3} 1) If for some neighbourhood $M \subseteq \Te$ of
$\tth$ we have both
(\ref{unifcon'}) and that for all $\ep>0$, $\PP\left( \inf\nolimits_{M, |\te - \tth| 
> \ep} |G(\te)|>0\right) = 1$, then any sequence
$\tten$ of $G_n$--estimators satisfies 
\bee
\label{unique}
\PP(\tten\in M, |\tten - \tth| > \ep)\to 0 \ \mbox{ for all }  \ep> 0.
\eee

\noindent
2)  Assume that  Condition \ref{conasy1} holds for all compact
subsets $M\subseteq\Te$, and that $\tth\in$ {\rm int}\,$\Te$ is the unique 
root of the equation $G(\te)=0$. Then the conditions in 1) hold for any compact
neighbourhood of $\tth$.

(a) Suppose $\Te$ is compact. Then any sequence of $G_n$--estimators is
weakly $\tth$--consistent, and this sequence is eventually unique. 

(b) Suppose $\Te$ is homeomorphic to $\R^p$, and choose any
$\rho\in\Te$. Then any sequence of $\rho$--centered
$G_n$--estimators is weakly $\tth$--consistent, and this sequence is
eventually unique, in the sense that if $\tten^\rho$ and $\tten'^{\rho'}$ 
are two sequences of, respectively, $\rho$ and $\rho'$--centered
$G_n$--estimators, then
$\PP(\tten^\rho\neq\tten'^{\rho'})\to1$. 
\end{thm}

Note that the asymptotic behavior of $\tten^\rho$ in (b) does 
not depend on the choice of $\rho$.

\subsection{Rate of convergence and asymptotic distribution} 

The results in the previous subsection ensure only the existence of a
solution converging to $\bar{\theta}$, but say nothing about the rate of the
convergence. To obtain results about the rate, we need a stronger condition
like, for instance, the following.

\begin{con}
\label{conasy2}
There exists a sequence of positive real numbers $a_n$ increasing
to infinity, such that the sequence of
random variables $a_n G_n(\bar{\theta})$ is stochastically bounded, i.e.\
such that for every $\ep > 0$ there exists a $K > 0 $ such that
$\PP \left( | a_n G_n(\bar{\theta}) | > K \right) < \ep$ for all $n$.
\end{con}

If the sequence $a_n G_n(\bar{\theta})$ converges in distribution, it is
stochastically bounded. The sequence $a_n$ can obviously be chosen in many
ways. As appears from the following theorem, the most interesting choice is
a sequence that goes to infinity as fast as possible.

\begin{thm}
\label{thmasy4}
Under Conditions \ref{conasy1} and \ref{conasy2} there is a sequence
$(\tten)$ of $G_n$--estimators such that
\bee
\label{rate2}
\lim_{c \rightarrow \infty}~ \sup_n
\PP \left( | a_n (\hat{\theta}_n - \bar{\theta} ) | \leq c \right) = 1
\eee
or in other words the sequence
$(a_n|\tten-\tth|)$ is stochastically bounded. Moreover, any sequence
$(\tten)$ of $G_n$--estimators which is weakly $\bar{\te}$--consistent satisfies
(\ref{rate2}).
\end{thm}

In some cases the coordinates of the estimator $\hat{\theta}_n$ do not
all converge to $\bar{\theta}$ at the same rate, see Section
\ref{highfreq}. When this happens, 
rates of convergence follow from the result below on the asymptotic
distribution of $G_n$--estimators, which can usually be derived from the
asymptotic distribution of the functions $G_n$ and their derivatives, as
stated in the next condition.

In that condition we have a $p$-dimensional variable $Z$, defined on an
extension of the space $(\Om,\f,\PP)$. By an extension we mean a
triple $(\bar{\Om},\bar{\f},\bar{\PP})$ with $\bar{\Om}=\Om\times\Om'$
and $\bar{\f}=\f\otimes\f'$ with $(\Om',\f')$ another measurable space,
and $\bar{\PP}$ is a probability measure satisfying $\bar{\PP}(A\times
\Om')=\PP(A)$ for all $A\in\f$. Any variable on $\Om$ or $\Om'$ is
extended in a trivial way as a variable on $\bar{\Om}$. Then we say
that a sequence $Z_n$ of $p$-dimensional random variables on $\Om$ converges
stably in law to $Z$ if $\E(Yh(Z_n))\rightarrow\bar{\E}(Yh(Z))$ for any
bounded variable $Y$ on $(\Om,\f)$ and any continuous bounded function $h$
on $\R^p$, and we write $Z_n\stackrel{\laa_{\rm st}}{\longrightarrow}Z$.
This automatically implies $Z_n\stackrel{\laa}{\longrightarrow}Z$.

\begin{con}
\label{conasymatrix}
There exist a sequence $A_n\in\m_p^{\rm inv}$ with each
entry of $A_n^{-1}$ tending to zero, a random vector $Z$ on an extension
of the space, and a random $\m_p$--valued function $W$, with
$W(\tth)$ almost surely in $\m_p^{\rm inv}$, such that 
for a neighbourhood $M$ of $\tth$ we have the following two properties:
\bee\label{C1}
A_nG_n(\tth)~\stackrel{\laa_{\rm st}}{\longrightarrow}~Z,
\eee
\bee\label{C2}
\sup_{\te\in M}\|\,A_n\partial_\te G_n (\te)A_n^{-1}-
W(\te)\| \stackrel{\PP}{ \rightarrow}~ 0.
\eee
\end{con}
In particular, this condition implies that
\bee\label{C3}
\left(\begin{array}{c}A_nG_n(\tth)\\
A_n\,\partial_\te G_n (\tth)\,A_n^{-1}\end{array}\right)~
\stackrel{\laa}{\longrightarrow}~
\left(\begin{array}{c} Z\\W(\tth)\end{array}\right).
\eee
When $W(\tth)$ is non-random, this is implied by (\ref{C2}), plus
$A_nG_n(\tth)~\stackrel{\laa}{\longrightarrow}~Z$ instead of (\ref{C1}).

\begin{thm}
\label{thmasy5}
Assume Conditions \ref{conasy1} and \ref{conasymatrix} holds, and let
$\tten$ be a weakly $\bar{\te}$--consistent sequence of $G_n$--estimators. Then
\bee
\label{asydis3}
A_n (\tten- \tth) \stackrel{{\cal L}_{st}}{\longrightarrow}
-W(\bar{\theta})^{-1} Z
\eee
and
\bee
\label{asydis3'}
A_n \partial_\theta G_n (\hat \theta_n) (\tten- \tth) \stackrel{{\cal
    L}_{st}} {\longrightarrow}
-Z.
\eee
\end{thm}

If the sequence $A_n^{-1}$ does not go fast enough to $0$, Condition
\ref{conasymatrix} may hold with $Z=0$, and (\ref{asydis3}) only gives a
rate of convergence which is not sharp. So, this result becomes
really interesting when Condition \ref{conasymatrix} holds with a
variable $Z$ that is {\em strongly non-degenerate}, in the sense that the
support of its law is not included in any proper linear subspace of
$\R^p$.

Quite often $Z$ is, conditionally on $\f$, centered Gaussian with an
invertible covariance matrix $V=V(\om)$. If $V$ is non-random this amounts
to having $Z$ independent of $\f$ and Gaussian. In the general
conditionally Gaussian case, the limit distribution in (\ref{asydis3})
is the normal variance-mixture with characteristic function 
$s \mapsto \E \left( \exp \left[ -\mbox{\small $\frac12$} s^*
W(\bar{\theta})^{-1} V W(\bar{\theta})^{*-1} s \right] \right),$ which
is Gaussian if both $V$ and $W(\bar{\theta})$ are non-random.
If one can construct
weakly consistent estimators $\widehat{V}_n$ for $V$, in the sense that
$\widehat{V_n}\stackrel{\PP}{\longrightarrow} V$ and each $\widehat{V}_n$
is positive definite, then (\ref{asydis3'}) implies that ($I_p$ is the
$p\times p$ identity matrix):
$$\widehat{V}_n^{-1/2}\,A_n\,\partial_\te G_n (\hat \theta_n)\,(\tten-\tth)
~\stackrel{\cal L} {\longrightarrow}~\n_p(0,I_p),$$
from which one can easily find confidence regions for the parameter
$\tth$. This may be difficult using the non-standard distribution in
(\ref{asydis3}).

\sect{Ergodic processes}
\label{asyergo}

In this subsection we consider the case where the observed process is
ergodic and present simple conditions that imply the previous general
assumptions. We assume that we 
have a sequence of random variables $X_1,X_2,\ldots$ with a
measurable state space $(D,\da)$, which is ergodic under the true
measure $\PP$. By ergodic we here mean that, for every integer $r\geq1$,
there is a probability measure $Q_r$ on $D^r$, such that for any function
$f:\, D^r \mapsto \R$ that is $Q_r$--integrable,
\bee
\label{ergo}
\frac1n \sum_{i=r}^n f(X_{i-r+1}, \ldots ,X_i)
\stackrel{\PP}{\longrightarrow} Q_r (f),
\eee
where $Q_r (f)$ denotes the integral of $f$ with respect to $Q_r$. It
is a weak form of ergodicity, which encompasses the case where we
observe a continuous time Markov process $Y$ at equidistant time points,
i.e. $X_n=Y_{n\De}$, $\De>0$. Suppose the state space of $Y$ is a domain
$D\subseteq \R^d$, and the transition kernels of $Y$ have
positive Lebesgue-densities $p_t(x,y)$, so $\PP(Y_t\in A\mid Y_0=x)
=\int_Ap_t(x,y)\,dy$ for all $t>0$ and $x\in D$. If $Y$ has a unique
invariant probability measure $\mu$, and if $p_t(x,y)$ is a continuous
function of $x$ for all $t>0$ and $y\in D$, then (\ref{ergo}) holds
for any initial measure $\eta$ on $D$. In this case, as an example,
$Q_2(dx,dy)=p_\De(x,y)\mu(dx)dy$.

We assume that at stage $n$ we observe the $n$ first variables
$X_1,\ldots,X_n$, and we consider estimating functions of the form
\bee
\label{form}
G_n(\theta) = \frac1n \sum_{i=r}^n g(X_{i-r+1}, \ldots ,X_i;\theta),
\eee
where $g: D^r \times \Theta \mapsto \R^p$ is jointly measurable and
satisfies the following assumption.

\begin{con}\label{conasy4} There is parameter value $\tth\in\mbox{\em
int}\,\Te$ and a neighbourhood $N$ of $\tth$ in $\Te$, such that:

\noindent
(1) The function $g(\theta): (x_1,\ldots,x_r) \mapsto g(x_1, \ldots,x_r;
\theta)$ is $Q_r$--integrable for all $\te\in N$, and
$Q_r \left( g(\bar{\theta}) \right) = 0$.

\noindent
(2) The function $\theta~ \mapsto~ g(x_1, \ldots,x_r; \theta)$ is $C^1$
on $N$ for all $(x_1,\ldots,x_r)$ in $D^r$.

\noindent
(3) For all compact subsets $M$ of $N$ there is a $Q_r$--integrable
function $\bar{g}_M$ on $D^r$ such that $|\partial_\te g(\te)|
\leq \bar{g}_M$ for all $\te\in M$.

\noindent
(4) The $p\times p$ matrix $Q_r(\partial_{\te}g(\tth))$ is invertible.
\end{con}
\vsq

Note that Condition \ref{conasy4} (3) is the property that
the function $\partial_\te g(\te)$ is locally dominated
$Q_r$--integrable for $\te \in N$. This is a traditional condition in
the statistical literature.

\begin{thm}
\label{thmasy6} 
Under Condition \ref{conasy4}, the estimating functions $G_n$ satisfy
Condition \ref{conasy1} for any compact subset $M \subseteq N$ and the
(non-random) $\R^p$-valued function $G(\te)=Q_r(g(\te))$. Hence there
exists an eventually unique sequence of weakly $\bar{\te}$--consistent 
$G_n$--estimators.

If further $N=\Te$ and $\tth$ is
the unique root of the function $G(\te)$, then any sequence $\tten$ of
$G_n$--estimators satisfies (\ref{unique}) for all compact
neighbourhoods $M$ of $\tth$. Moreover, 
if $\Te$ is compact, then any sequence $\tten$ of $G_n$--estimators is
weakly $\tth$-consistent, and if $\Te$ is homeomorphic to $\R^p$, then
for any $\rho\in\Te$, any  sequence $\tten^\rho$  of $\rho$-centered
$G_n$--estimators is weakly $\tth$-consistent.
\end{thm}

Suppose that we assumed ergodicity in a slightly stronger sense, namely
that in (\ref{ergo}) the convergence takes place $\PP$--almost surely.
This is, for instance, the case for a discretely
observed continuous time Markov process that satisfies the conditions
given above. Then in Condition \ref{conasy1} we obtain almost sure
convergence as well, and by Theorem \ref{thmasy2} we see that there exists
a sequence of strongly $\bar{\te}$--consistent $G_n$--estimators.

Finally, let us assume that the estimating functions $G_n$ also
satisfy a central limit theorem:
\bee\label{erg1}
\sqrt{n}\, G_n(\bar{\theta}) = \frac{1}{\sqrt{n}}\sum_{i=r}^n g(X_{i-r+1},
\ldots ,X_i;\bar{\theta}) \stackrel{\cal L}{\longrightarrow} N(0,
V(\bar{\theta}))
\eee
for some (necessarily non--random) $p \times p$--matrix $V(\tth)$. Then it
follows from Theorem \ref{thmasy5} that any sequence $\tten$ of
weakly $\tth$--consistent $G_n$--estimators satisfies
\[\sqrt{n}\,(\tten - \tth)~ \stackrel{\cal L}{\longrightarrow}~
\n_p\big(0,Q_r(\partial_\te g(\tth))^{-1}\, V(\tth)\,
Q_r(\partial_\te g(\tth))^{*-1}\big).
\]

\sect{Longitudinal data}
\label{repeat}
 
The previous section applies, in particular, when the observations $X_i$
are i.i.d., the situation for which estimating functions were
initially introduced. In this case one naturally takes $r=1$, (\ref{ergo})
is the usual LLN (with a.s.\ convergence) and (\ref{erg1}) the usual CLT.
However, each observation $X_i$ may have a complex structure. An
interesting example is longitudinal data, where each $X_i$ consists of
observations of a stochastic process. This kind of data, that are also
referred to as panel data, have received a lot of attention in the
statistical and econometric literature. In this branch of statistics,
estimating functions have frequently been applied. A classical text is
\cite{diggle}. For applications of diffusion processes in a
longitudinal data context, see \cite{Pedersen2000} and \cite{ditlevsen}.

As an example, we consider the case where the $X_i$'s are $n$ independent 
copies of a discretely observed stationary Markov process $Y$. More 
specifically, we have a state space $D$, and for each $\te$ a stationary
Markov process $Y^\te$ with transition semi-group $(P^\te_t)_{t\geq0}$ and 
stationary initial distribution $\mu^\te$.  The $i$th longitudinal
observation is $X_i=(Y^i_{j\De}:\,j=0,1, \ldots,m)$, where $\De>0$,
$m$ is a fixed integer, and $Y^1,\ldots,Y^n$ are i.i.d.\ copies of
$Y^{\te_0}$, where $\te_0$ is the true parameter value. 

We consider the estimating function
$$G_n(\te)=\frac1n\sum_{i=1}^n\sum_{j=1}^mh(Y^i_{(j-1)\De},Y^i_{j\De};\te),$$
where $h$ is a measurable function on $D^2\times\Te$ which satisfies
\bee\label{D5}
\int_Dh(x_1,x_2;\te)\,P^\te_\De(x_1,dx_2)=0
\eee
for all $\te\in\Te$ and all $x_1 \in D$. This means that the inner sum
defines a {\em martingale estimating function} for the $i$th Markov
process. In contrast to the previous section, $m$ is here fixed, and
we exploit the martingale property only to compute the limiting
variance. Below, $Q(dx_1,dx_2)=\mu^{\te_0}(dx_1)\,P^{\te_0}(x_1,dx_2)$ is 
the law of $(Y^1_0,Y^1_\De)$ under the true measure. Besides (\ref{D5}),
we impose the following condition on $h$:

\begin{con}
\label{cc1}
$ \ $ \\
\noindent
(1) The function $h(\te):~(x_1,x_2)\mapsto h(x_1,x_2;\te)$ is
square--integrable with respect to $Q$ for all $\te\in\Te$.

\noindent
(2) The function $\te\mapsto h(x_1,x_2; \te)$ is $C^1$ for all $(x_1,x_2) 
\in D^2$.

\noindent
(3) For all compact subsets $M$ of $\Theta$ there is a $Q$--integrable
function $\bar{h}_M$ on $D^2$ such that $|\partial_\te h(\te)|
\leq \bar{h}_M$ for all $\te\in M$.

\noindent
(4) The $p \times p$ matrix $Q(\partial_\te h(\te_0))$ is invertible.
\end{con}

This seems identical to Condition \ref{conasy4} with $r=2$ and
$\tth=\te_0$ ($Q(h(\te_0))=0$ follows from (\ref{D5})), but the meaning 
of $Q$ here is quite different from the meaning of $Q_2$ in that condition.

\begin{thm}
\label{thmasy8}
Under Condition \ref{cc1} there exists a sequence $\tten$ of strongly
consistent $G_n$--estimators, eventually unique and
satisfying (as $n\rightarrow\infty$)
\bee\label{D4}
\sqrt{n}\,(\tten-\te_0)~\stackrel{\laa}{\longrightarrow}~
\n_p\big(0,m^{-1}\,V(\te_0)\big),
\eee
where
$V(\te_0)=\big(Q(\partial_\te h(\te_0))\big)^{-1}\,
Q(h(\te_0)\,h(\te_0)^*)\,\big(Q(\partial_\te h(\te_0))\big)^{\star-1}.$

If further $Q\big(g(\te)\big)\neq0$ for all $\te\neq\te_0$, then any
sequence $\tten$ of $G_n$--estimators satisfies (\ref{unique}) with
$\bar \te = \te_0$ for all compact
neighbourhoods $K$ of $\te_0$. Moreover,  if $\Te$ is compact, then any
sequence $\tten$ of $G_n$--estimators satisfies (\ref{D4}), and
if $\Te$ is homeomorphic to $\R^p$, then for any $\rho\in\Te$ and
any sequence $\tten^\rho$ of $\rho$--centered $G_n$--estimators satisfies 
(\ref{D4}).
\end{thm}

\sect{High frequency observation of a diffusion}
\label{highfreq}

In this section we present examples which illustrate that it is
necessary to allow the limiting functions in Conditions \ref{conasy1}
and \ref{conasymatrix} to be random, and that it is necessary to allow
different rates of convergence for the coordinates of the estimators. 

For each value $\te=(\al,\be)$ in a subset $\Te=A\times B$ of $\R^2$, we 
consider the one--dimensional diffusion process given by the stochastic 
differential equation
\bee\label{H1}
dX^\te_t=a(X^\te_t,\al)\,dt+b(X^\te_t;\be)\,dW_t,
\qquad X^\te_0=x_0,
\eee
with $W$ a standard Wiener process on a filtered probability space
$(\Om,\f,(\f_t)_{t\geq0},\PP)$, where $\f$ is the $\PP$--completion of
$\bigvee_{t>0}\f_t$ and $(\f_t)$ the filtration generated by $W$. 

We observe $X=X^{\te_0}$ at the times $i\De_n,\,i=0,1,\ldots,n$, for
the true parameter value $\te_0=(\al_0,\be_0)$, and we assume that
$\De_n\to0$ as $n\to\infty$. Thus we have observation in the time interval
$[0,T_n]$ only, where the total time span is $T_n=n\De_n$. Since $X_0$ is
observed, it is no restriction to assume that the starting point is a
non-random number $x_0$, independent of $\te$.

We make the following smoothness assumptions on the coefficients. These 
could be substantially weakened, at the price of more complex proofs. That 
the state space is the entire set $\R$ could also be relaxed.

\begin{assu}\label{H2}
$ \ $ \\
\noindent(1) The function $b$ is $C^3$ on $\R\times B$, and all its
derivatives are of polynomial growth in $x$, uniformly in $\be\in K$ for all
compact subsets $K\subseteq B$.

\noindent(2) The function $a$ is $C^2$ on $\R\times A$, and all its
derivatives are of polynomial growth in $x$, uniformly in $\al\in K$ for all
compact subsets $K\subseteq A$.

\noindent(3) The functions $a(\cdot;\al)$ and $b(\cdot;\beta)$ are
globally Lipschitz for all $\te=(\al,\be)$.

\noindent(4) ~$\inf_{x\in\R,\,\be\in K}\,b(x;\be)>0$ for all compact
subsets $K \subseteq B$.
\end{assu}

This implies in particular that (\ref{H1}) has a unique strong
solution, which is Markov, and also that for any $m\geq0$, $t \geq 0$
and any compact subset $K\subseteq\Te$:
\bee\label{H3}
\sup_{\te\in K,\,0\leq\De\leq1}\,\De^{-m/2}\,
\E\big(\sup\nolimits_{s\in[0,\De]}\,|X^\te_{t+s}-X^\te_t|^m
\mid\f_t\big)<\infty.
\eee

The situation is quite different when $T_n\equiv T$ for some fixed $T>0$
and when $T_n\to\infty$ as $n\to\infty$. Below, we consider
these two cases.

\subsection{Fixed time span}

When $T_n\equiv T$, it is well known that the drift cannot be estimated
consistently, so we consider only estimation of the second component
$\be$ of the parameter $\te$. In other words, the set $A$ consists of
a single point, or $a(x,\al)=a(x)$ does not depend on a parameter. The
drift coefficient still satisfies the relevant conditions in Assumption \ref{H2},
but it may be unknown (the estimating functions (\ref{H4}) does not
depend on $a$), so we are solving a semi-parametric problem.

With the notation $c(x,\be)=b(x,\be)^2$ and
$\De^n_iX=X_{i\De_n}-X_{(i-1)\De_n}$, a simple estimating function for
$\be$, which yields the estimator proposed by \cite{genon93}, is given by
\bee\label{H4}
G_n(\be)=\sum_{i=1}^n\frac{\partial_\be c(X_{(i-1)\De_n};\be)}
{c(X_{(i-1)\De_n};\be)^2}\Big((\De^n_iX)^2-\De_n\,c(X_{(i-1)\De_n};\be)\Big).
\eee

\begin{thm}\label{TD} Suppose $T_n\equiv T$ and that Assumption \ref{H2} 
and the identifiability condition $\int_0^T|\partial_\te b(X_t,\be_0)|
\,dt>0$ outside a $\PP$-null set hold. Then the estimating function $G_n$ 
satisfies the Conditions \ref{conasy1} and \ref{conasymatrix} with $\tth=\te_0
=\be_0$, any compact subset $M\subseteq B$ and the random functions
\bee\label{H5}
G(\be)=\int_0^T\frac{\partial_\be c(X_t;\be)}{c(X_t;\be)^2}\big(
c(X_t;\be_0)-c(X_t;\be)\big)\,dt,
\eee
\bee\label{H6}
W(\be)=\partial_\be G(\be)=
\int_0^T\left[\partial_\be\Big(\frac{\partial_\be c(X_t;\be)}
{c(X_t;\be)^2}\Big)\big(c(X_t;\be_0)-c(X_t;\be)\big)
-\frac{(\partial_\be c(X_t;\be))^2}{c(X_t;\be)^2}\right]\,dt.
\eee
Moreover, $A_n=\rn$ and $Z$ is a random variable which
conditionally on $\f$ is centered Gaussian with variance $-2T\,W(\be_0)$.

Suppose further that $G(\be)\neq0$ for all $\be\neq\be_0$, then any
sequence $\widehat{\be}_n$ of $G_n$--estimators satisfies (\ref{unique}) with
$\tten = \widehat{\be}_n$ and $\bar \te = \be_0$ for all compact
neighbourhoods $K$ of 
$\be_0$. Moreover, if $B$ is compact, then any sequence $\widehat{\be}_n$ 
of  $G_n$--estimators satisfies (\ref{asydis3}) and (\ref{asydis3'}), and
if $B$ is an open (finite or infinite) interval, then for all $\rho\in
B$, any sequence $\widehat{\be}_n^\rho$ of $\rho$--centered
$G_n$--estimators satisfies (\ref{asydis3}) and (\ref{asydis3'}). 
\end{thm}

The estimator studied here is efficient and has the optimal rate of
convergence, see \cite{Gobet01}.  A general theory for
approximate martingale estimating functions when $T_n\equiv T$ can be
found in \cite{jakobsen2017}. Another general class of 
estimators was investigated in \cite{genon93}.

\subsection{Time span going to infinity}

Now, we assume that $T_n\to\infty$, so drift parameters can be
consistently estimated, and we consider the full parameter space $\Te=A\times B$.
We need the following ergodicity assumption.

\begin{assu}\label{H2E} Under the true parameter value $\te_0$ the solution
of (\ref{H1}) is ergodic, in the sense that there is a probability measure
$\mu_{\te_0}$ on $\R$ (necessarily an invariant measure for the Markov process
$X^{\te_0}$) such that any $\mu_{\te_0}$--integrable function $f$ satisfies
\bee\label{HH3}
\frac1n\sum_{i=0}^nf(X^{\te_0}_{i\De_n})\stackrel{\PP}{\longrightarrow}
\mu_{\te_0}(f).
\eee
Moreover, $\sup_{t>0}\,\E(|X^{\te_0}_t|^m)<\infty$ for all $m\geq0$, and hence
$\int|x|^m\,\mu_{\te_0}(dx)<\infty$ as well.
\end{assu}

Let us also state an identifiability assumption.

\begin{assu}\label{H2I'} We have
$\int_\R|\partial_\be b(x,\be_0)|\,\mu_{\te_0}(dx)>0$ and
$\int_\R|\partial_\al a(x,\al_0)|\,\mu_{\te_0}(dx)>0$.
\end{assu}

We use the two-dimensional estimation function $G_n=(G_n^1,G_n^2)$ given 
for $\te=(\al,\be)$ by 
\bee\label{H6_2}
\begin{array}{l}
G_n^1(\te)=\frac1{T_n}\sum\limits_{i=1}^n
\frac{\partial_\al a(X_{(i-1)\De_n};\al)}
{c(X_{(i-1)\De_n};\be)}\Big(\De^n_iX-\De_n\,a(X_{(i-1)\De_n};\al)\Big)\\
G_n^2(\te)=\frac1{T_n}\sum\limits_{i=1}^n\frac{\partial_\be 
c(X_{(i-1)\De_n};\be)}
{c(X_{(i-1)\De_n};\be)^2}\Big((\De^n_iX)^2-\De_n\,c(X_{(i-1)\De_n};\be)\Big).
\end{array}
\eee
Note that $G_n^2$ equals the estimating
functions (\ref{H4}).

\begin{thm}\label{TE} Suppose $T_n\rightarrow\infty$ and $n\De_n^2
\rightarrow0$, and that Assumptions \ref{H2}, \ref{H2E} and \ref{H2I'} hold.
Then the estimating function $G_n$ given by (\ref{H6_2}) satisfies
Conditions \ref{conasy1} and \ref{conasymatrix} for any compact subset
$M\subseteq \Te$ and with $\tth=\te_0$, $A_n=\mbox{\rm
  diag}(\sqrt{T_n},\rn)$, and the non-random functions 
\bean
G^1(\te)&=&\int_\R \frac{\partial_\al a(x;\al)}{c(x;\be)}\big(
a(x;\al_0)-a(x;\al)\big)\,\mu_{\te_0}(dx), \\
G^2(\te)&=&\int_\R \frac{\partial_\be c(x;\be)}{c(x;\be)^2}\big(
c(x;\be_0)-c(x;\be)\big)\,\mu_{\te_0}(dx) 
\eean
\bean
W^{11}(\te)&=&\partial _\alpha G^1 (\te) =
\int_\R\left[\partial_\al\Big(\frac{\partial_\al a(x;\al)} 
{c(x;\be)}\Big)\big(a(x;\al_0)-a(x;\al)\big)
-\frac{(\partial_\al a(x;\al))^2}{c(x;\be)}\right]\,\mu_{\te_0}(dx) \\
W^{22}(\te)&=& \partial _\be G^2 (\te) = \int_\R\left[ \partial_\be
\Big(\frac{\partial_\be c(x;\be)} {c(x;\be)^2}\Big)\big(c(x;\be_0)-c(x;\be)\big)
-\frac{(\partial_\be c(x;\be))^2}{c(x;\be)^2}\right]\,\mu_{\te_0}(dx) \\ 
W^{12}(\te)&=& \partial _\be G^1 (\te) = -\int_\R\frac{\partial_\al
  a(x;\al)\,\partial_\be c(x;\be)}
{c(x;\be)^2}\,\big(a(x;\al_0)-a(x;\al)\big)\,\mu_{\te_0}(dx)  \\
W^{21}(\te)&=& \partial _\alpha G^2 (\te) = 0,
\eean
and with $Z=(Z^1,Z^2)$ a  two-dimensional variable independent of
$\f$, where $Z^1$ and $Z^2$ are independent centered Gaussian with
variances $-W^{11}(\te_0)$ and  $- W^{22}(\te_0)$.

Suppose further that $G(\te)\neq0$ for all $\te\neq\te_0$, then any
sequence $\tten$ of $G_n$--estimators satisfies (\ref{unique}) with
$\bar \te = \te_0$ for all compact neighbourhoods $K$ of
$\te_0$. Moreover, if $\Te$ is compact, then any sequence $\tten$ of
$G_n$--estimators satisfies (\ref{asydis3}) and (\ref{asydis3'}), and
if $\Te$ is homeomorphic to $\R^2$, then for all $\rho\in\Te$,
any sequence  $\tten^\rho$ of $\rho$--centered $G_n$--estimators 
satisfies (\ref{asydis3}) and (\ref{asydis3'}).
\end{thm}

The estimator studied here is efficient and has an optimal rate of
convergence, see \cite{Gobet02}. A general theory for approximate
martingale estimating functions for diffusion processes under the
asymptotic scenario considered in this subsection can be found in
\cite{efficient}. 

\sect{Proofs}
\label{proofs}

The proof of Theorem \ref{thmasy1} is based on the fixed point theorem
for contractions. A mapping $f$ from a subset $M$ of $\R^p$ into $\R^p$ is called a {\it contraction} if there exists a constant $C\in(0,1)$, called
a contraction constant, such that $|f(x)-f(y)|\leq C|x-y|$ for all
$x,y\in M$. A proof of the following lemma can for instance be found
at page 229 in \cite{LoomisSternberg}.  Below, $\bar{B}_r(x)$ is the
closed ball of $\R^p$ with radius $r$ and center $x$.

\begin{lem}
\label{lemasy2}
Let $f: \bar{B}_r(x_0) \mapsto \R^p$ be a contraction such that
$ | f(x_0) - x_0 | \leq (1 - C)r$, where $C$ is the contraction
constant. Then $f$ has a unique fixed point $x$ (i.e., $f(x)=x$) in
$\bar{B}_r(x_0)$, which is the limit of the sequence $x_n$ started
at $x_0$ and defined by induction through $x_{n+1}=f(x_n)$.
\end{lem}

Specifically, the fixed point theorem is used to prove the following lemma.

\begin{lem}
\label{lemasy1}
Let $f$ be a differentiable mapping from a closed subset $M$
of $\R^p$ into $\R^p$, and let $A\in\m_p^{\rm inv}$.
Define $\lambda=\frac12 \| A^{-1} \|^{-1}$. If
\[
\| \partial_x f(x) - A \|\leq\lambda
\]
on a ball $\bar{B}_r(x_0)\subseteq M$, then any point $y$ in $\bar{B}_{\la
r}(f(x_0))$ is the image $y=f(x)$ of a unique point $x$ in $\bar{B}_r(x_0)$.
\end{lem}

\pr  Choose $y \in \bar{B}_{\lambda r}(f(x_0))$ and define the function
$\phi (x) = x + A^{-1} (y - f(x))$.  It is sufficient to prove that
$\phi$ has a unique fixed point in $\bar{B}_r(x_0)$. Since
$\partial_x\phi(x)=I-A^{-1}\partial_x f(x)$
for all $x \in \bar{B}_r(x_0)$, it follows that
\[
\|\partial_x\phi(x)\|\leq\|A^{-1}\|\|A-\partial_xf(x)\|\leq
\mbox{\small $\frac12$}.
\]
Thus $\phi$ is a contraction on $\bar{B}_r(x_0)$ with contraction constant
$1/2$. Since
\[
| \phi(x_0) - x_0 | = | A^{-1} (y - f(x_0)) | \leq
\| A^{-1} \| \lambda r = r/2,
\]
the result follows from the previous lemma. \qed
\vsq

\noindent
{\bf Proof of Theorem \ref{thmasy1}}. 1) Let us introduce the random variables
$$ \Lambda=\mbox{\small $\frac12$}\|\partial_\te G(\tth)^{-1}\|^{-1},$$
$$Y(\ep)=\sup_{\te:|\te-\tth|\leq\ep}~\|\partial_\te G(\te)
-\partial_\te G(\tth)\|,$$
$$Z_n=\sup_{\te\in M}\|\partial_\te G_n(\te)-\partial_\te G(\te)\|,$$
which are $\f$--measurable (because $\f$ is $\PP$-complete), and
also the $\f$--measurable sets
$$C_{n,\ep}=\left\{Y(\ep)\leq{\mbox{\small $\frac12$}\Lambda}\right\}\cap
\left\{Z_n\leq{\mbox{\small $\frac12$} \Lambda}\right\}\cap
\left\{|G_n(\bar{\te})|\leq\Lambda\ep\right\}.$$ 

On the set $C_{n,\ep}$ we have $Y(\ep)+Z_n\leq\La$, hence
$\|\partial_\te G_n(\te)-\partial_\te G(\tth)\|\leq\La$ whenever
$|\te-\bar{\te}|\leq\ep$, and also $|G_n(\tth)|\leq\La\ep$. Thus for
any given $\om$ in $C_{n,\ep}$
we can apply Lemma \ref{lemasy1} with $f=G_n$ and $\la=\Lambda$ and
$A=\partial_\te G(\tth)$ and $r=\ep$, to get
\bee\label{4exis}
\om\in C_{n,\ep}~~\Rightarrow~~\mbox{there is a unique~
$\te_{n,\ep}(\om)\in \bar{B}_\ep(\bar{\te})$~ with}
~~G_n(\te_{n,\ep}(\om),\om)=0.
\eee
Moreover by Lemma \ref{lemasy2} and the proof of Lemma \ref{lemasy1},
for each $\om\in C_{n,\om}$ we have $\te_{n,\ep}=\lim_p\,z_p$ for the
sequence defined inductively by $z_0=\bar{\te}$ and $z_{p+1}=z_p-
\partial_\te G(\tth)G_n(z_p)$. Hence $\te_{n,\om}$ (restricted to
$C_{n,\ep}$) is $\f_n$--measurable.

2) Let us now prove the existence of a sequence $\ep_n\downarrow0$ which
satisfies
\bee\label{4lim}
\PP(C_{n,\ep_n})~\to~1.
\eee
For this, we first recall the well known fact that a sequence
of real--valued variables $V_n$ goes to $0$ in probability if and only if
there is a sequence $\ep_n\downarrow0$ such that $\PP(|V_n|\geq\ep_n)\to0$.

Condition  \ref{conasy1} yields $\PP(\Lambda>0)=1$ and
$G_n(\bar{\te})\stackrel{\PP}{ \rightarrow}~ 0$ and
$Z_n\stackrel{\PP}{\rightarrow}~ 0$ and $\limsup_{\ep\to0}\,Y(\ep)=0$
almost surely. Hence $G_n(\bar{\te})/\La\stackrel{\PP}{ \rightarrow}~ 0$
and $Z_n/\La\stackrel{\PP}{\rightarrow}~ 0$ and $\limsup_{\ep\to0}\,Y(\ep)/\La=0$ almost surely, and we deduce that, for some sequence
$\ep_n\downarrow0$,
\bee\label{4lim2}
\PP(C_{n,\ep_n})\geq
1-\PP\left(Y(\ep_n)/\La>{\mbox{\small$\frac12$}}\right)-
\PP\left(Z_n/\La>{\mbox{\small $\frac12$}}\right)
-\PP\left(|G_n(\tth)|/\La>\ep_n\right)~\rightarrow~1,
\eee
which yields (\ref{4lim}).

3) Now we are ready to prove the existence of a weakly
$\tth$--consistent sequence $\tten$ of $G_n$--estimators. First, we
choose arbitrary $G_n$--estimators $\tten'$, which are
known to exist. Then, with $\ep_n$ as above, we 
define $\tten$ to be equal to $\te_{n,\ep_n}$ on the set $C_{n,\ep_n}$
and to $\tten'$ on the complement of this set. This gives us an
$\f_n$-measurable variable $\tten$, whose weak $\tth$--consistency
readily follows from the fact that $|\te_{n,\ep_n}-\tth|\leq\ep_n$, plus
$\ep_n\rightarrow0$ and (\ref{4lim}).

4) It remains to prove the last claim. We assume that we have two
sequences $(\tten)$ and $(\tten')$ of $G_n$--estimators, both of them weakly 
$\bar{\te}$-consistent, and we want to prove $\PP(\tten\neq\tten')\to0$. Since
 $\tten-\tth$ and $\tten'-\tth$ go to $0$ in probability, we can find a
 sequence $\ep_n\downarrow0$ such that (\ref{4lim}) holds, together with
\bee\label{4lim1}
\PP(|\tten-\tth|\geq\ep_n)~\to~0,\qquad
\PP(|\tten'-\tth|\geq\ep_n)~\to~0.
\eee

Since both $\tten$ and $\tten'$ solve the estimating equation
$G_n(\te)=0$ when this equation has a solution, we readily deduce from
(\ref{4exis}) that on the set $C_{n,\ep_n}$ and if
$|\tten-\tth|\leq\ep_n$ and $|\tten'-\tth|\leq\ep_n$, then necessarily
$\tten=\tten'$. Hence $\PP(\tten\neq\tten')\to0$ follows from (\ref{4lim})
and (\ref{4lim1}).

\qed
\vsc

\noindent
{\bf Proof of Theorem \ref{thmasy3}}. 
By $d$ we denote the Euclidean distance, except in case 2 (b)  where
$d$ denotes the distance used in (\ref{CC1}). In both cases
$d(\te,\de)=1$ for all $\te\in\Te$ by convention. For any  
neighbourhood $M\subseteq \Te$ of $\tth$ and $\ep\in(0,1)$, we define 
$Z(M)_n=\sup_{\theta \in M}\,|G_n(\te)-G(\te)|$  and 
$Y(M,\ep)=\inf_{\te\in M,\,d(\te,\tth)>\ep}\,|G(\te)|$. Under the
assumptions in 1),  $Y(M,\ep)>0$ and
$Z(M)_n\stackrel{\PP}{\longrightarrow}0$, so as $n\to\infty$:
\bee\label{CC2}
\PP\big(Z(M)_n < Y(M,\ep)\big)~\rightarrow~1.
\eee 
Because $\{ \tten \in M \} \subseteq \{ Z(M)_n \geq |G(\tten)| \}$, then
any sequence $\tten$ of $G_n$--estimators satisfies 
\bee\label{CC3}
d(\tten,\tth)\leq\ep~~\mbox{on the set}~~
\{\tten\in M\}\cap\{Z(M)_n < Y(M,\ep)\}.
\eee 
This implies the conclusion in 1).

The conditions in 2) obviously implies that the conditions in 1) hold for any
compact neighbourhood of $\tth$.
In case (a), (\ref{CC2}) and (\ref{CC3}) with $M=\Te$, plus the facts that 
$\{\tten\in\Te\}=D_n$ and $\PP(D_n)\to1$ by Theorem \ref{thmasy1}, yield 
the weak $\tth$-consistency of $\tten$. In case (b), choose a weakly 
$\tth$--consistent sequence $\tten$ of $G_n$--estimators (which exists
by Theorem \ref{thmasy1}), thus $\PP(B_{\ep,n})\to1$, where $B_{\ep,n}=D_n\cap\{
d(\tten,\tth)\leq\ep\}$. Moreover $d(\tten^\rho,\rho)\leq d(\tten,\rho)$ by 
(\ref{CC1}), hence $\tten^\rho\in M_\ep$ on the set $B_{\ep,n}$, where
$M_\ep$ is the compact set $\{\te\in\Te:\,d(\te,\rho)\leq 
d(\tth,\rho)+\ep\}$. Then (\ref{CC2}) and (\ref{CC3}) with $M=M_\ep$ and
$\tten^\rho$ imply $\PP(d(\tten^\rho,\tth)>\ep)\to0$, hence the weak 
$\tth$--consistency of $\tten^\rho$. Finally, in both cases, the
eventual uniqueness follows from Theorem \ref{thmasy1}. \qed
\vsc

\noindent
{\bf Proof of Theorem \ref{thmasy4}}.
Since each variable $a_n|\tten-\tth|$ is finite--valued (recall the
convention that $|\tten-\tth|=1$ when $\tten=\de$), it is well known that
(\ref{rate2}) is in fact equivalent to the apparently weaker
requirement that
\bee
\label{rate2'}
\lim_{c \rightarrow \infty}~ \limsup_n
\PP \left( | a_n (\hat{\theta}_n - \bar{\theta} ) | \leq c \right) = 1.
\eee

We take for $\tten$ the estimator constructed in Theorem
\ref{thmasy1}, and in the following we use the notation of the proof of this
theorem. Condition \ref{conasy2} yields that for any $\ep>0$ there is a
$K>0$ such that $\PP(|G_n(\tth)|>K\La/a_n)\leq\ep$ for all $n$. Then
(\ref{4lim2}) and the fact that $Y(\eta)\to0$ as $\eta\to0$ and that
$Z_n\to0$ as $n\to\infty$, in probability, yield that
\bee\label{4lim4}
\liminf_n\PP(C_{n,K/a_n})\geq1-\ep.
\eee
Moreover, with the sequence $\ep_n$ constructed in part 3 of the proof
of Theorem \ref{thmasy1}, we know that on $C_{n,\ep_n}$ we have
$|\tten-\tth|\leq\ep_n$. Hence by (\ref{4exis}) we necessarily have
$|\tten-\tth|\leq K/a_n$ as well on the intersection $C_{n,\ep_n}\cap
C_{n,K/a_n}$. Then, combining (\ref{4lim}) and (\ref{4lim4}), we
deduce that
$$\liminf_n \PP(|\tten-\tth|\leq K/a_n)\geq1-\ep$$
and, since $\ep>0$ is arbitrary, (\ref{rate2'}) readily follows.

Finally, the last claim follows from the eventual uniqueness
proved in Theorem \ref{thmasy1}.
\qed
\vsc

\noindent
{\bf Proof of Theorem \ref{thmasy5}}.
By shrinking $M$ if necessary, we may  suppose that $M$ is convex. Define
$C_n=\{\tten\in M\}$. By the weak $\tth$-consistency we have
$\PP(C_n)\to1$, and the mean value theorem yields that on $C_n$
\[
G_n(\tten)-G_n(\tth)=\partial_\te \widetilde G_n (\tten-\tth).
\]
Here $\partial_\te \widetilde G_n$ is the $p\times p$--matrix whose
$jk$th entry is $\partial_\te G_n(\te^{(j)}_n)_{jk}$, where each
$\te^{(j)}_n$ is a (random) convex combination of $\tten$ and
$\tth$. Observe that
\[
\|A_n\,\partial_\te \widetilde G_n \,A_n^{-1} -W(\tth)\|
\leq p \left(\sup_{\te\in M}\|\,A_n\,\partial_\te G_n (\te)\,A_n^{-1}
-W(\te)\| + \hspace{-2mm} \sup_{\te:\,|\te-\tth| \leq |\tten
-\tth|} \hspace{-1mm} \|W(\te)-W(\tth)\| \right)
\]
on $C_n$. Hence $A_n\,\partial_\te \widetilde G_n A_n^{-1}
\stackrel{\PP}{\longrightarrow} W(\tth)$ because $W$ is a.s.\ continuous
and $\tten$ is weakly $\tth$-consistent. Therefore, on the set
$C_n \cap \{A_n\,\partial_\te \widetilde G_n A_n^{-1}$ is invertible$\}$ 
(the probability of which goes to $1$) we have $G_n(\tten)=0$ and thus
$$\begin{array}{lll}
A_n(\tten-\tth)&=&-(A_n\,\partial_\te \widetilde G_n A_n^{-1})^{-1}\,
A_nG_n(\tth) \\ &=&-W(\tth)^{-1}\,A_nG_n(\tth)+U_n A_nG_n(\tth),
\end{array}$$
where $U_n\stackrel{\PP}{\longrightarrow}0$. In view of (\ref{C1}), this
yields (\ref{asydis3}), and since the convergence is stable and 
$A_n\,\partial_\te G_n(\tten)\,A_n^{-1} \stackrel{\PP}{\longrightarrow}W(\tth)$, 
(\ref{asydis3'}) also follows. \qed
\vsc

\noindent
{\bf Proof of Theorem \ref{thmasy6}}.
Define $G(\te)=Q_r(g(\te))$. Our hypotheses on $g$, the dominated
convergence theorem (implying in particular $\partial_\te G(\te))=
Q_r(\partial_\te g)$ when $\te\in N$) and (\ref{ergo}) clearly yield
all requirements of Condition \ref{conasy1}, except (\ref{unifcon}).

For proving (\ref{unifcon}), we first deduce from (\ref{ergo}) and
$\partial_\te G(\te)=Q_r(\partial_\te g)$ that, for any $\te\in N$,
\bee\label{D1}
\partial_\te G_n(\te)~\stackrel{\PP}{\longrightarrow}~
\partial_\te G(\te).
\eee
Next, we define for $\eta>0$ a function $k_\eta$ on $D^r$ and a real number
$\al_\eta$ by
\bean
&&k_\eta(x_1, \ldots , x_r) =
\sup_{\te,\te'\in M: |\te'-\te| \leq \eta}
\| \partial_{\theta} g(x_1, \ldots , x_r; \te') -
\partial_{\theta} g(x_1, \ldots , x_r;\theta) \|\\
&&\al_\eta=\sup_{\te,\te'\in M: |\te'-\te| \leq \eta}
\| \partial_{\te} G(\te') -\partial_{\te} G(\te) \|.
\eean
By (2) of Condition \ref{conasy4} the functions $\te\mapsto \partial_\te
g(s_1,\ldots,x_r;\te)$ are uniformly continuous on the compact set $M$,
hence (3) of this Condition and the dominated convergence theorem yield
$Q_r(k_\eta)\to0$ as $\eta\to0$. The function $\partial_\te G(\te)$ is
continuous on the compact set $M$, so $\al_\eta\rightarrow0$.

By the finite covering property of the compact set $M$, for any
$\eta>0$, we have a partition of $l(\eta)$ nonempty subsets $M_j$ of $M$
with diameters less than $\eta$. For each $j=1,\ldots,l(\eta)$ we choose a
point $\te_j\in M_j$, and we set
$$B_\eta^n=\sum_{j=1}^{l(\eta)}\|\partial_\te G_n(\te_j)-\partial_\te
G(\te_j)\|,
\qquad Z^n_\eta=\frac1n\sum_{i=r}^nk_\eta(X_{i-r+1},\ldots,X_i).$$

Recalling that $\partial_\te G_n(\te)$ equals
$\frac1n\sum_{i=r}^n\partial_\te g(X_{i-r+1},\ldots,X_i;\te)$, we
see that
$$\te\in M_j~~\Rightarrow~~
\|\partial_\te G_n(\te)-\partial_\te G(\te)\|\leq
\|\partial_\te G_n(\te_j)-\partial_\te G(\te_j)\|+Z^n_\eta+\al_\eta,$$
hence
\bee\label{D2}
\sup_{\te\in M}\,\|\partial_\te G_n(\te)-\partial_\te G(\te)\|\leq
B_\eta^n+Z^n_\eta+\al_\eta,
\eee
By (\ref{ergo}) and (\ref{D1}) we have $Z^n_\eta\stackrel{\PP}
{\longrightarrow}Q_r(k_\eta)$ and $B^n_\eta\stackrel{\PP}
{\longrightarrow}0$, both as $n\rightarrow\infty$ and for any fixed $\eta$.
Thus, for any $\ep>0$ we can choose first $\eta>0$ and secondly $n_0$ large
enough to have
$$\al_\eta\leq\ep,\qquad Q_r(k_\eta)\leq\ep,\qquad n\geq n_0~~
\Rightarrow~~\PP(B^n_\eta>\ep)+\PP(Z^n_\eta>Q_r(k_\eta)+\ep)\leq\ep.$$
From this and (\ref{D2}), we readily deduce the following, which gives us (\ref{unifcon}):
$$\PP\big(\sup_{\te\in M}\,\|\partial_\te G_n(\te)-\partial_\te G(\te)\|
>4\ep\big)~\leq~\ep.$$

The second part of the theorem follows readily from Theorem \ref{thmasy3}. 

\qed
\vsc

\noindent
{\bf Proof of Theorem \ref{thmasy8}}.
Let $\widetilde{Q}$ denote the law of $(Y^{\te_0}_0,\ldots,Y^{\te_0}_{m\De})$ on
$D^{m+1}$. With
$$g(x;\te)=\sum_{j=1}^mh(y_{(j-1)},x_{j};\te)\quad\mbox{for}~~
x=(y_0,\ldots,y_m)\in D^{m+1}$$
we have $G_n(\te)=\frac1n\sum_{i=1}^ng(X_i;\te)$. Our hypotheses
readily imply that $g$ satisfies Condition \ref{conasy4} with $N=\Te$,
$\tth =\theta_0$, $r=1$ and $Q_1=\widetilde{Q}$. The $X_i$'s are
i.i.d.\ and hence satisfy (\ref{ergo}) with almost sure
convergence. Thus the existence of strongly consistent, eventually
unique $G_n$--estimators and the last part of the theorem follow from
Theorem \ref{thmasy6} and the comments which follow it. Note that
(\ref{unifcon}) for a compact neighbourhood $M$ of $\te_0$ also
follows, with $G(\te)=\widetilde{Q}(g(\te))=mQ(h(\te))$. 

The central limit theorem for i.i.d.\ variables yields (\ref{C1}) with
$\tth=\te_0$ and $A_n=\sqrt{n}\,I_p$ ($I_p$ is the identity matrix), 
where $Z$ is independent of $\f$
and centered Gaussian with covariance matrix $V(\te_0)=
\widetilde{Q}(g(X_{1};\te_0)g(X_{1};\te_0^*))$, which by
(\ref{D5}) is equal to $mQ(h(\te_0)h(\te_0)^*)$. Moreover 
(\ref{unifcon}) and $A_n=\sqrt{n}\,I_p$ imply (\ref{C2}) with $W(\te)=
\partial_\te G(\te)=mQ(\partial_\te h(\te))$. Now (\ref{D4}) follows
from Theorem \ref{thmasy5}, and the last statements follow from
Theorem \ref{thmasy3}. \qed

\vsc

\noindent
{\bf Proof of Theorem \ref{TD}}. 1) Below, $M\subseteq B$ is a
fixed compact subset, and $C$ is a generic constant. In view of
Assumption \ref{H2} and of (\ref{H3}), the following estimates are
classical (and easy to derive: use It\^o's formula applied to $X^2W$
and to $c(X;\be_0)$ for the third one). 
\bean
&&\big|\,\E((\De^n_iX)^2\mid\f_{(i-1)\De_n})-\De_nc(X_{(i-1)\De_n};\be_0)\big|
\leq C\De_n^2\big(1+|X_{(i-1)\De_n}|^C\big) \\
&&\E(|\De^n_iX|^{2m}\mid\f_{(i-1)\De_n})\leq C_m\De_n^m
\big(1+|X_{(i-1)\De_n}|^{C_m}\big) \rule{0mm}{5mm}\\
&&\big|\,\E((\De^n_iX)^2\,\De^n_iW\mid\f_{(i-1)\De_n})\big|
\leq C\De_n^2\big(1+|X_{(i-1)\De_n}|^C\big). \rule{0mm}{5mm}
\eean
Then, setting
\bean
f(x,\be)&=&\frac{\partial_\be c(x;\be)}
{c(x;\be)^2}\big(c(x;\be_0)-c(x;\be)\big)\\
\ze(\be)^n_i&=&\frac{\partial_\be c(X_{(i-1)\De_n};\be)}
{c(X_{(i-1)\De_n};\be)^2}\big((\De_i^nX)^2-\De_nc(X_{(i-1)\De_n};\be)\big),
\eean
($\ze(\be)^n_i$ is the $i$th summand in (\ref{H4})), we deduce
for all $\be\in M$ and $m=0,1$:
\bea
&&\Big|\,\E(\partial^m_\be\ze(\be)^n_{i}\mid\f_{(i-1)\De_n})
-\De_n\partial^m_\be f(X_{(i-1)\De_n};\be)
\Big|\leq C\De_n^2\big(1+|X_{(i-1)\De_n}|^{C}\big)~~  \label{1}\\
&&\E(|\partial^m_\be\ze(\be)^n_{i}|^2\mid\f_{(i-1)\De_n})\leq C\De_n^2
\big(1+|X_{(i-1)\De_n}|^{C}\big)\label{2} \rule{0mm}{5mm}\\
&&\E\big(\sup\nolimits_{\be\in M}\,|\partial^2_\be\ze(\be)^n_{i}|
\mid\f_{(i-1)\De_n}\big)\leq C\De_n\big(1+|X_{(i-1)\De_n}|^{C}\big)
\label{3} \rule{0mm}{5mm} \\
&&\Big|\,\E((\ze(\be_0)^n_{i})^2\mid\f_{(i-1)\De_n})-2\De_n^2
\frac{\big(\partial_\be c(X_{(i-1)\De_n};\be_0)\big)^2}
{c(X_{(i-1)\De_n};\be_0)^2}\Big|\leq C\De_n^3
\big(1+|X_{(i-1)\De_n}|^C\big)\qquad\label{4} \rule{0mm}{5mm} \\
&&\E(|\ze(\be_0)^n_{i}|^4\mid\f_{(i-1)\De_n})\leq C\De_n^4
\big(1+|X_{(i-1)\De_n}|^C\big)\label{5} \rule{0mm}{5mm} \\
&&\Big|\,\E(\ze(\be_0)^n_{i}\De^n_{i}W\mid\f_{(i-1)\De_n})\Big|
\leq C\De_n^2\big(1+|X_{(i-1)\De_n}|^C\big) \rule{0mm}{5mm} \label{6}
\eea

2) Define $G$ and $W$ by 
(\ref{H5}) and (\ref{H6}). Then $G(\be_0)=0$, the $C^1$ property of $G_n$
and $G$ are obvious, $W(\be)=\partial_\be G(\be)$, and the
identifiability condition implies $W(\be_0)<0$ a.s. We have $\De_n=T/n$, hence
(\ref{1}) and (\ref{2}) for $m=0,1$ plus Riemann integration yield $G_n(\be_0)
\stackrel{\PP}{\longrightarrow}0$ and $\partial_\be G_n(\be)
\stackrel{\PP}{\longrightarrow}W(\be)$. Inequality (\ref{3}) implies that 
$\be\mapsto\partial_\be G_n(\be)$) is Lipschitz on $M$ with a (random) 
Lipschitz coefficient that is integrable. This together with the
pointwise in $\be$ convergence in probability implies by 
standard tightness arguments that (\ref{unifcon}) and (\ref{C2}) with
$A_n=\rn$ hold.

By Theorem IX.7.28 of \cite{JacodShiryaev} plus again Riemann
integration, (\ref{C1}) is a straightforward consequence of (\ref{1})
with $m=0$ and $\be=\be_0$, plus (\ref{4})--(\ref{6}). Hence, we have
proved that Conditions \ref{conasy1} and \ref{conasymatrix} hold, and the
last claim follows from Theorem \ref{thmasy3}. \qed
\vsc

\noindent
{\bf Proof of Theorem \ref{TE}}. 
Let us denote by $\ze'^j(\te)^n_i$ the $i$th summands in the
definition (\ref{H6_2}) of $G^j_n(\te)$, $j=1,2$. We also define for $\te=(\al,\be)$
\[
\widetilde f(x,\te)=\frac{\partial_\al a(x;\al)}{c(x,\be)}\,
\left(a(x;\al_0)-a(x;\al) \right)
\]
and observe that
$$\big|\,\E(\De^n_iX\mid\f_{(i-1)\De_n})-\De_na(X_{(i-1)\De_n};\al)\big|
\leq C\De_n^2\big(1+|X_{(i-1)\De_n}|^C\big).$$
Then for all $\te\in M$ (with $M$ a compact subset of $\Te$) and $m=0,1$
we have
\bean
&&\hskip-1cm\Big|\,\E(\partial^m_\te\ze'^1(\te)^n_{i}\mid\f_{(i-1)\De_n})
-\frac{\De_n}{T_n}\,\partial^m_\te \widetilde f(X_{(i-1)\De_n};\te)
\Big|\leq C\frac{\De_n^2}{T_n}\,\big(1+|X_{(i-1)\De_n}|^C\big)\\
&&\hskip-1cm\E(|\partial^m_\te\ze'^1(\te)^n_{i}|^2\mid\f_{(i-1)\De_n})\leq
C\frac{\De_n}{T_n^2}\,\big(1+|X_{(i-1)\De_n}|^C\big)\\
&&\hskip-1cm\E\big(\sup\nolimits_{\te\in M}\,|\partial^2_\te
\ze'^1(\te)^n_{i}|\mid\f_{(i-1)\De_n}\big)\leq C\frac{\De_n}{T_n}\,
\big(1+|X_{(i-1)\De_n}|^C\big)\\
&&\hskip-1cm\Big|\,\E((\ze'^2(\te_0)^n_{i})^2\mid\f_{(i-1)\De_n})-\frac{\De_n}{T_n^2}
\,\frac{(\partial_\al a(X_{(i-1)\De_n};\al))^2}
{c(X_{(i-1)\De_n};\be)}\Big|\leq
C\frac{\De_n^2}{T_n^2}\,\big(1+|X_{(i-1)\De_n}|^C\big)\\ 
&&\hskip-1cm\E(|\ze'^1(\te_0)^n_{i}|^4\mid\f_{(i-1)\De_n})\leq
C\frac{\De_n^2}{T_n^4}\,\big(1+|X_{(i-1)\De_n}|^{C}\big)\\
&&\hskip-1cm\big|\,\E((\ze'^1(\te_0)^n_i\ze'^2(\te_0)^n_i\mid\f_{(i-1)
\De_n})\big|\leq C\frac{\De_n^2}{T_n^2}\,\big(1+|X_{(i-1)\De_n}|^C\big).
\eean
Since $\ze'^2(\te)^n_i=\ze(\be)^n_i/T_n$, the variables $\ze'^1(\te)^n_i$
satisfy (\ref{1})--(\ref{6}), after normalization by an appropriate
power of $T_n$. 

At this stage, the proof follows the scheme of the previous proof, except
that we use the law of large numbers (\ref{HH3}) instead of the convergence
of Riemann sums (we leave the -- tedious -- details to the reader). In 
particular, (\ref{C1}) is a consequence of classical 
convergence results, for which we  do not need the analogue of (\ref{6}) 
because in the ergodic case the limit $Z$ in (\ref{C1}) is 
automatically independent of $\f$. Hence,
Conditions \ref{conasy1} and \ref{conasymatrix} hold, and the
last claim follows from Theorem \ref{thmasy3}.\qed

\bibliographystyle{natbib}
\bibliography{asystat}

\end{document}